\input amstex

\documentstyle{amsppt}
\document
\topmatter
\title
  QCH  K\"ahler surfaces II
\endtitle
\author
W\l odzimierz Jelonek
\endauthor

\abstract{ In this paper we give new examples of QCH K\"ahler surfaces whose opposite almost Hermitian structure is Hermitian and not locally conformally K\"ahler.  In this way we give also a large class of examples of Hermitian surfaces with $J$-invariant Ricci tensor which are not l.c.k.}
 \endabstract
\endtopmatter
\define\r{\rightarrow}
\define\G{\Gamma}
\define\DE{\Cal D^{\perp}}

\define\n{\nabla}
\define\om{\omega}
\define\w{\wedge}
\define\k{\diamondsuit}
\define\th{\theta}
\define\p{\partial}
\define\a{\alpha}

\define\lb{\lambda}

\define\1{D_{\lb}}
\define\2{D_{\mu}}
\define\bO{\bar{\Omega}}

\define\0{\Omega}
\define\Om{\Omega}

\define\bJ{\bar J}

\define\J{\Cal J}
\define\wg{\wedge}
\define\De{\Cal D}
\define\R{\Cal R}
\define\dl{\delta}



{\bf 0. Introduction.}  In our paper [J-2] we have studied QCH semi-symmetric K\"ahler surfaces $(M,g,\bJ)$ which are ambi-K\"ahler,  i.e.,  which admit an opposite locally conformally K\"ahler Hermitian structure. We gave a classification of such K\"ahler semi-symmetric surfaces.  The natural question was whether there exist semi-symmetric K\"ahler surfaces whose opposite almost Hermitian structure is Hermitian and not l.c.k.
The aim of the present paper is to give new examples of semi-symmetric K\"ahler surfaces $(M,g,\bJ)$  foliated by Euclidean spaces. Our examples admit a Hermitian opposite structure $J$ which satisfies the condition  $G_1$ of Gray and which is {\bf not} locally conformally K\"ahler.  Hermitian surface  $(M,g,J)$ has $J$-invariant Ricci tensor,  the conformal scalar curvature $\kappa$ of $(M,g,J)$ is equal to the scalar curvature $\tau$ of $(M,g)$, hence the zero defect $\tau-\kappa$ of $(M, g,J)$ is $0$.  For these surfaces also  $\rho^*=\rho$, where $\rho^*$ is the *-Ricci tensor and $\rho$ is the Ricci tensor of $(M,g,J)$.  In the last part of the paper we also give examples of QCH K\"ahler surface (see [J-1]) which are not semi-symmetric and whose opposite almost Hermitian structure is Hermitian and not locally conformally K\"ahler.  In this way we give a large class of Hermitian surfaces with $J$-invariant Ricci tensor such that $J$ is not l.c.k. In the paper we use some results from [J-4]  which we cite here with proofs for convenience of the reader.
\par
\bigskip
{\bf 1.  Hermitian 4-manifolds. } Let $(M,g,J)$ be an almost Hermitian manifold, i.e., an  almost complex structure $J$ is  orthogonal with respect to $g$: $g(X,Y)=g(JX,JY)$ for all $X,Y\in\frak X(M)$. We say that $(M,g,J)$ is a Hermitian  manifold if its almost Hermitian structure $J$ is integrable, which means that the Niejenhuse tensor $N^J(X,Y)$ vanishes. This is also equivalent to integrability of a complex distribution $T^{(1,0)}M\subset TM\otimes\Bbb C$.
Let us write $\wg^2M=\wg^2T^*M$. In the sequel we shall identify the bundle $TM$ with $T^*M$ by means of $g$, so we shall also write $\wg^2M=\wg^2TM$.
The Hodge star operator $*$ (which depends on the orientation of $M$) defines an endomorphism $*:\w ^2M\rightarrow\w^2 M$ with $*^2=id$ and we denote by $\w^+,\w^-$ its eigensubbundles corresponding to $1,-1$ respectively.
In the sequel we shall consider 4-dimensional Hermitian manifolds $(M,g,J)$, which we shall also call Hermitian surfaces. Such manifolds are always oriented and we choose an orientation in such a way that the K\"ahler form $\0(X,Y)=g(JX,Y) $ is a self-dual form (i.e., $\0\in\w^+M$). The vector bundle of self-dual forms admits a decomposition
$$\w^+M=\Bbb{R}\0\oplus LM,\tag 1.1$$
where by $LM$ we denote the bundle of real $J$-skew invariant 2-forms (i.e., $LM=\{\Phi\in\w M:\Phi(JX,JY)=-\Phi(X,Y)$\}). The bundle $LM$ is a complex line bundle over $M$ with the complex structure $\Cal J$ defined by $(\Cal J\Phi)(X,Y)=-\Phi(JX,Y)$.
For a Hermitian surface the covariant derivative of the K\"ahler form $\0$ is locally expressed by
$$\n \0=a\otimes\Phi+\J a\otimes\J\Phi, \tag 1.2$$
where $\J a(X)=-a(JX)$. The Lee form $\th$ of $(M,g,J)$ is defined by the equality $d\0=2\th\w \0$. We have $2\th=-\delta\0\circ J$. By $\rho$ we denote the Ricci tensor of a Riemannian manifold $(M,g)$ and by $\tau$ the scalar curvature of $(M,g)$, i.e.,  $\tau=\text{ tr}_g\rho$.
 A Hermitian manifold $(M,g,J)$ is said to have $J$-invariant Ricci tensor $\rho$ if $\rho(X,Y)=\rho(JX,JY)$ for all $X,Y\in \frak X(M)$.   An involutive distribution is called foliation. A foliation $\De$ is called totally geodesic if its every leaf is a totally geodesic submanifold of $(M,g)$, i.e., $\n_XY\in \G(\De)$ for every   $X,Y\in\G(\De)$.
 A Hermitian 4-manifold $(M,g,J)$ is said to have an opposite Hermitian structure  if it admits an orthogonal Hermitian structure $\bar J$ with an anti-self-dual K\"ahler form $\bO$.
For any almost Hermitian four manifold the following formula holds (see [G-H]) :
$$\frac12(\rho(X,Y)+\rho(JX,JY))-\frac12(\rho^*(X,Y)+\rho^*(Y,X))=\frac14(\tau-\tau^*)g(X,Y), \tag 1.3$$
where $\rho^*$ is the Ricci *-tensor defined by
$$\rho^*(X,Y)=\frac12\text{tr }\{Z\rightarrow R(X,JY)JZ\}.\tag 1.4$$
 The curvature tensor is $R(X,Y)Z=([\n_X,\n_Y]-\n_{[X,Y]})Z$.    A Hermitian manifold has a $J$-invariant Ricci tensor if and only if
  $$\rho(X,Y)-\rho^*(X,Y)=\frac{\tau-\kappa}6g(X,Y)$$  (see [M]).  The function $\tau-\kappa$ is called the zero defect of  $(M,g,J)$.
  By $\De$ we denote the nullity distribution of $(M,g,J)$ defined by $\De=\{X\in TM: \n_XJ=0\}$. For a Hermitian manifold it follows from (1.2) that $\De$ is $J$-invariant. Consequently  dim$\De=2$ in $M$ if $|\n J|\ne0$ in $M$. In the sequel we shall assume that $\n J\ne0$ on the whole of $M$.  By $\DE$ we shall denote the orthogonal complement of $\De$ in $M$.  An almost Hermitian structure $\bJ$ defined by $\bJ X=JX$ if $X\in \DE$  and $\bJ X=-JX$ if $X\in\De$ is called a natural opposite almost Hermitian structure.  In the sequel we shall study Hermitian surfaces $(M,g,J)$  with $\n J\ne0$ on the whole of $M$, whose natural opposite structure $\bJ$ is K\"ahler.

The curvature tensor $R$ of a 4-dimensional manifold $(M,g)$ determines an endomorphism $\R$ of the bundle $\w^2 M$ defined by $g(\R(X \w Y),Z\w W)=\R(X\w Y,Z\w W)=- R(X,Y,Z,W)=-g(R(X,Y)Z,W)$. The conformal scalar curvature $\kappa$ is defined by  $$\kappa=\tau-6(|\th|^2+\dl\th).\tag 1.6$$
We say that an almost Hermitian manifold $(M,g,J)$ satisfies the second condition $(G_2)$ of A. Gray if its curvature tensor $R$ satisfies the condition
$$R(X,Y,Z,W)-R(JX,JY,Z,W)=R(JX,Y,JZ,W)+R(JX,Y,Z,JW),\tag $G_2$ $$
for all $X,Y,Z,W\in \frak X(M)$. We say that it satisfies the condition $(G_3)$ of A. Gray if
$$R(JX,JY,JZ,JW)=R(X,Y,Z,W),\tag $G_3$ $$
 for all $X,Y,Z,W\in \frak X(M)$.
Let us define $ B=\frac12(\R-*\R*); W=\frac12(\R+*\R*)_{0}=\frac12(\R+*\R*)-\frac{\tau}{12}Id; W^+=\frac12(W+*W);W^-=\frac12(W-*W)$. Then
$$\R=\frac{\tau}{12}Id+B+W^++W^-.$$
The tensor $W$ is called the Weyl tensor and its components $W^+,W^-$ are called the self-dual and anti-self-dual Weyl tensors.
The  QCH  K\"ahler surface is a K\"ahler surface $(M,g,\bJ)$ with distribution $\De$ such that the holomorphic curvature
$K(\pi)=R(X,JX,JX,X)$ of any $J$-invariant $ 2$-plane $\pi\subset
T_xM$, where $X\in \pi$ and $g(X,X)=1$, depends only on the point
$x$ and the number $|X_{\De}|=\sqrt{g(X_{\De},X_{\De})}$, where
$X_{\De}$ is the orthogonal projection of $X$ onto $\De$. In this
case  we have
$$R(X,JX,JX,X)=\phi(x,|X_{\De}|),$$ where $\phi(x,t)=a(x)+b(x)t^2+c(x)t^4$ and
 $a,b,c$ are smooth functions on $M$. Also $R=a\Pi+b\Phi+c\Psi$
 for certain curvature tensors $\Pi,\Phi,\Psi\in \bigotimes^4\frak X^*(M)$
  of K\"ahler type.  Every QCH K\"ahler surface admit an opposite almost Hermitian structure $J$ which satisfies the second Gray condition.   A QCH surface is semi-symmetric if  $J$ satisfies the first Gray condition
$$R(JX,JY,Z,W)=R(X,Y,Z,W)\tag $G_1$ $$
or equivalently if $\tau-\kappa=0$, where $\kappa$  is a conformal curvature of $(M,g,J)$,  since $R.R=\frac{\tau-\kappa}6\Pi.R$  (see [J-1]).
A foliation $\Cal F$
on a Riemannian manifold $(M,g)$ is called conformal if for every
$V\in\G(T\Cal F)$
$$L_Vg=\a(V)g$$ holds on $T\Cal F^{\perp}$, where $\a$ is a one
form vanishing on $T\Cal F^{\perp}$ (see [J-3]).
\bigskip
{\bf 2. Hermitian surfaces with Hermitian Ricci tensor.} In the first part of this paragraph we give lemmas A,B,C, D with proofs, which were first proved in [J-4].
 Note that
for every manifold satisfying condition $(G_3)$ we have $\R (LM)\subset \w^+M$, its Ricci tensor $\rho$ is $J$-invariant and its $*$-Ricci tensor is symmetric. Indeed, since $R(j(X\w Y),j(Z\w W))=R(X\w Y,Z\w W)$, where $j(X\w Y)=JX\w JY$, we have $\R(\text{ ker }(j-id),\text{ ker }(j+id))=0$. Since ker$(j-id)=\w^-M\oplus\Bbb{R}\0$ and ker$(j+id)=LM$, we get $g(\R(LM),\w^-M\oplus\Bbb{R}\0)=0$.  Consequently $\R(LM)\subset LM\subset\w^+M$. In fact, the condition $\R(LM)\subset\w^+M$ holds if and only if the Ricci tensor $\rho$ of $(M,g)$ is $J$-invariant ( see [D] p.5 (i)) and an almost Hermitian 4-manifold $(M,g,J)$ with $J$-invariant Ricci tensor and symmetric *-Ricci tensor satisfies $(G_3)$.
\par
\bigskip
{\bf Lemma A . }{\it Let $(M,g,J)$ be a Hermitian 4-manifold. Let us assume that  $|\n J| \ne 0$ on $M$. Then for any local orthonormal oriented basis $\{E_1,E_2\}$ of $\DE$ there exists a global oriented orthonormal basis $\{E_3,E_4\}$ of $\De$ independent of the choice of $\{E_1,E_2\}$  such that
$$\n \0=\a(\th_1\otimes\Phi+\th_2\otimes\Psi), \tag 2.1$$
where  $\Phi=\th_1\w\th_3-\th_2\w\th_4,\ \Psi=\th_1\w\th_4+\th_2\w\th_3$,$\a=\frac1{2\sqrt{2}}|\n J|$ and $\{\th_1,\th_2,\th_3,\th_4\}$ is a cobasis dual to $\{E_1,E_2,E_3,E_4\}$.  Moreover $\delta\0=-2\a\th_3,\th=-\a\th_4$.}
\medskip
{\it Proof. } Let $\{E_1,E_2\}$ be any orthonormal basis of $\DE$, $E_2=JE_1$ . Then (1.2) holds, where $a=\a\th_1$. Let us choose any orthonormal basis $\{E_3',E_4'=JE_3'\}$ in $\De$.
Let us define $\Phi'=\th_1\w\th_3'-\th_2\w\th_4',\ \Psi'=\th_1\w\th_4'+\th_2\w\th_3'$. Then
$\{\Phi',\Psi'\}$ is an oriented orthonormal local basis in $LM$. Thus we have
$$\Phi=\cos\phi \Phi'-\sin\phi\Psi',\ \ \Psi=\sin\phi\Phi'+\cos\phi\Psi'.$$
for some local function $\phi$. Then $\n \0=\a\{\th_1(\cos\phi \Phi'-\sin\phi\Psi')+\th_2(\sin\phi\Phi'+\cos\phi\Psi')\}$. Let us define $E_3=\cos\phi E_3'-\sin\phi E_4',\ E_4=\sin\phi E_3'+\cos\phi E_4'$. Then $\{E_3,E_4\}$ is the basis we were searching for.
From (2.1) it is easy to get:$\delta\0=-2\a\th_3,\th=-\a\th_4$.$\k$
\medskip
Any frame $\{E_1,E_2,E_3,E_4\}$ constructed as above we shall call standard ( or special).
\medskip
The following Lemma is well known ( it means that for a Hermitian surface the component $W_3^+$ of the positive Weyl tensor vanishes).
\bigskip
{\bf Lemma B .} {\it Let $(M,g,J)$ be a Hermitian surface. Then for any local orthonormal basis $\{\Phi,\Psi\}$ of LM we have $\R(\Phi,\Phi)=\R(\Psi,\Psi)$ and $\R(\Phi,\Psi)=0$.}
\bigskip
 It is known  that a Hermitian manifold $(M,g,J)$ satisfies the second condition of Gray if and only if its Ricci tensor is $J$-invariant, it has symmetric $*$-Ricci tensor  and the component $W_3^+$ of positive Weyl tensor vanishes (i.e., $\R_{LM}=a id_{LM}$, where $\R_{LM}=p_{LM}\circ \R_{|LM}$ and $p_{LM}$ is the orthogonal projection $p_{LM}:\w M\rightarrow LM$). It is well known that any almost  Hermitian manifold satisfying condition $(G_2)$ satisfies the condition $(G_3)$ and that any Hermitian manifold satisfying $(G_3)$ satisfies $(G_2)$ (i.e., for Hermitian manifolds these two conditions are equivalent).

\par
\bigskip
{\bf Lemma C. }{\it Let $(M,g,J)$ be a Hermitian surface with $J$-invariant Ricci tensor (i.e., $\R(LM)\subset \w^+M$). Let $\{E_1,E_2,E_3,E_4\}$ be a local orthonormal frame such that (2.1) holds. Then
$$\gather\G^3_{11}=\G^3_{22}= E_3\ln \a,  \tag a\\
\G^3_{44}=\G^4_{21}= -\G^4_{12}=- E_3\ln \a, \tag b\\
\G^3_{21}=-\G^3_{12},\ \G^4_{11}=\G^4_{22},  \tag c\\
-\G^3_{21}+\G^4_{22}=\a,\tag d\\
\G^4_{33}=-E_4\ln\a+\a, \tag e\endgather$$
where $\n_XE_i=\sum\om^j_i(X)E_j$ and $\G^i_{kj}=\om^i_j(E_k)$.}

{\it Proof.}  Note that $\G^i_{kj}=-\G^j_{ki}$. We have $$\gather g(\n_{E_1}JX,Y)=\a\Phi(X,Y), g(\n_{E_2}JX,Y)=\a\Psi(X,Y)\tag 2.2\\
\n_{E_3}J=0,\n_{E_4}J=0.\endgather$$
Let us write $p(X)=\frac12g(\n_X\Phi,\Psi)=\om^2_1(X)+\om^4_3(X)$. Then we have
$ \n_X\0=\a\th_1(X)\Phi+\a\th_2(X)\Psi,
\n_X\Phi=-\a\th_1(X)\0+p(X)\Psi,
\n_X\Psi=-\a\th_2(X)\0-p(X)\Phi$.
Consequently using (2.2) we get $$\gather g(R(E_1,E_3).JX,Y)=-\n_{[E_1,E_3]}\0-E_3\a\Phi-\a p(E_3)\Psi,\tag 2.3a\\
g(R(E_1,E_4).JX,Y)=-\n_{[E_1,E_4]}\0-E_4\a\Phi-\a p(E_4)\Psi,\tag 2.3b\\
g(R(E_2,E_3).JX,Y)=-\n_{[E_2,E_3]}\0-E_3\a\Psi+\a p(E_3)\Phi,\tag 2.3c\\
g(R(E_4,E_2).JX,Y)=-\n_{[E_4,E_2]}\0+E_4\a\Psi-\a p(E_4)\Phi,\tag 2.3d\\
g(R(E_1,E_2).JX,Y)=-\n_{[E_1,E_2]}\0+(E_1\a-\a p(E_2))\Psi-(\a p(E_1)+E_2\a)\Phi,\tag 2.3e\endgather$$
where as usual $R(X,Y).J=\n_X(\n_YJ)-\n_Y(\n_XJ)-\n_{[X,Y]}J$. Recall that  $R(X,Y).J=R(X,Y)\circ J-J\circ R(X,Y)$, i.e., $R(X,Y)$ acts on  the tensor $J$ as a derivation.
 Since $\R(LM)\subset \w^+M$, it is clear that
$$\gather
g(R(E_1,E_3).JX,Y)=g(R(E_4,E_2).JX,Y),\tag 2.4a\\
g(R(E_3,E_2).JX,Y)=g(R(E_4,E_1).JX,Y).\tag 2.4b
\endgather$$
Consequently we get from (2.3) and (2.4) using the condition $\R(LM)\subset \w^+M$
$$\gather \frac12\R(\Phi,\Psi)=-g(R(E_1,E_3).JE_1,E_3)=(E_3\a+\a\th_1([E_1,E_3]))=E_3\a-\a\G^3_{11},\tag 2.5a\\
\frac12\R(\Phi,\Psi)=-g(R(E_2,E_3).JE_3,E_2)=-(E_3\a+\a\th_2([E_2,E_3]))=-E_3\a+\a\G^3_{22},\tag 2.5b\\
\frac12\R(\Phi,\Phi)=-g(R(E_4,E_2).JE_3,E_2)=(E_4\a-\a\th_2([E_4,E_2]))=E_4\a-\a\G^4_{22},\tag 2.5c\\
\frac12\R(\Psi,\Psi)=-g(R(E_1,E_4).JE_1,E_3)=-(-E_4\a-\a\th_1([E_1,E_4]))=E_4\a-\a\G^4_{11}.\tag 2.5d\\
\frac12\R(\Phi,\Psi)=g(R(E_1,E_4).JE_1,E_4)=-(\a p(E_4)+\a\th_2([E_1,E_4]))=\a\G^3_{44}-\a\G^2_{14},\tag 2.5e\\
\frac12\R(\Phi,\Psi)=-g(R(E_4,E_2).JE_1,E_3)=\a p(E_4)+\a\th_1([E_4,E_2]))=-\a\G^1_{24}-\a\G^3_{44},\tag 2.5f\\
\frac12\R(\Phi,\Phi)=g(R(E_1,E_3).JE_1,E_4)=-\a\th_2([E_1,E_3])-\a p(E_3)=-\a\G^2_{13}-\a\G^4_{33},\tag 2.5g\\
\frac12\R(\Psi,\Psi)=-g(R(E_2,E_3).JE_1,E_3)=\a\th_1([E_2,E_3])-\a p(E_3)=\a\G^1_{23}-\a\G^4_{33}.\tag 2.5h
\endgather$$ Note that for a Lee form $\th$ of $(M,g,J)$ we have $\th=-\a\th_4$. Let us write $X=2\a E_3$.
Then $d\0=-2\a\th_4\w\0$ and
$$L_X\0=d(i_X\0)+i_X(d\0)=2d(\a\th_4)=-2d\th.$$
Hence
$$L_X\0^2=2L_X\0\w\0=-4d\th\w\0.$$
Since $d\0=2\th\w\0$, we have $d\th\w\0=0$. Thus $L_X\0^2=0$ and div$X=0$. It means that $$\G^3_{11}+\G^3_{22}+\G^3_{44}=E_3\ln\a.\tag 2.6$$
From Lemma B, (2.5) and (2.6) we get equations (a),(b),(c) of Lemma C.
Since $2\th=-\dl\0\circ J$, we have $\n J(E_1,E_1)=\n J(E_2,E_2)=\a E_3;\n J(E_3,E_3)=\n J(E_4,E_4)=0$. Consequently we get equations
$$\a E_3+J(\n_{E_1}E_1)=\n_{E_1}E_2;\ \a E_3+J(\n_{E_2}E_2)=-\n_{E_2}E_1.$$
Thus $\G^3_{12}+\G^4_{11}=\a$ and $-\G^3_{21}+\G^4_{22}=\a$ and (d) follows. On the other hand
$\R(\Phi,\Phi)+\R(\Psi,\Psi)=4(E_4\a-\a\G^4_{11})=4(\a\G^1_{23}-\a\G^4_{33})$ and (e) follows.
$\k$

The first part of the next Lemma is well known (see [A-G-1]).

\par
\bigskip
{\bf Lemma D. }{\it Let $(M,g,J)$ be a Hermitian surface   with  J-invariant Ricci tensor. Then $\G^4_{13}=-E_2\ln\a, \ \ \G^4_{23}=E_1\ln\a$ and $d\th$ is anti-self-dual. In particular if $M$ is compact, then $(M,g,J)$ is locally conformally K\"ahler. In addition $(M,g,J)$ is l.c.K. if and only if $E_3\a=0,\ \G^1_{34}=E_2\ln\a, \ \G^2_{34}=-E_1\ln\a$.}
\medskip
{\it Proof. } From [A-G-1] Th.2 it follows that $\0$ is an eigenform of $W^+$. Let $\{E_1,E_2,E_3,E_4\}$ be the local special frame. From (2.3e) we obtain that an equation $\R(\0,\Phi)=0 $ is equivalent to  $-\a\th_2([E_1,E_2])-\a p(E_2)+E_1\a=0$. Analogously an equation $\R(\0,\Psi)=0 $ is equivalent to $-\a\th_1([E_1,E_2])-(\a p(E_1)+E_2\a)=0$. We get  $\G^4_{13}=-E_2\ln\a, \ \ \G^4_{23}=E_1\ln\a$  after some easy computation.  Note that
$$\gather d\th_4(E_3,E_4)=-\th_4([E_3,E_4])=E_3\ln\a,\tag 2.7a\\
d\th_4(E_1,E_2)=-\th_4([E_1,E_2])=-\G^4_{12}+\G^4_{21}=-2E_3\ln\a,\tag 2.7b\\
d\th_4(E_1,E_3)=-\th_4([E_1,E_3])=-\G^4_{13}+\G^4_{31}, \tag2.7c\\
d\th_4(E_2,E_3)=-\th_4([E_2,E_3])=-\G^4_{23}+\G^4_{32}.\tag 2.7d
\endgather$$
We also have  $d\th=-d\a\w\th_4-\a d\th_4$.
From (2.7) we get
$$-d\th=-2E_3\a\bO+(-\a\G^1_{34}+E_2\a)\bar{\Phi}+(-\a\G^4_{32}+E_1\a)\bar{\Psi},\tag 2.8$$
where $\bar{\Phi}=\th_1\w\th_3+\th_2\w\th_4,\bar{\Psi}=\th_1\w\th_4-\th_2\w\th_3$.
Consequently $d\th$ is anti-self-dual.$\k$

\medskip
If $(M,g,J)$ is a Hermitian surface with $|\n J|\ne 0$ on $M$, then the distributions $\De,\DE$ define a natural opposite almost Hermitian structure $\bar{J}$ on $M$. This structure is defined
as follows $\bar J|_{\De}=-J|_{\De},\bar J|_{\DE}=J|_{\DE}$. In the special basis we just have: $\bar J E_1=E_2, \bar J  E_3=-E_4$.

\par
\bigskip

{\bf Lemma E. } {\it Let $(M,g,J)$ be a Hermitian  4-manifold with Hermitian Ricci tensor  and K\"ahler natural opposite structure. Let us assume that $|\n J|\ne 0$ on $M$. Then

(a) $\De$ is a totally geodesic foliation,

(b)$E_3\ln\a=0$

(c) $\n_{E_4}E_4=0$.}

 \medskip
{\it Proof. }   Let us  choose a local orthonormal frame $\{E_1,E_2,..,E_4\}$ such that (2.1) holds.
Note that (we write $\n_X\th_i=\om_i^j(X)\th_j$, $\bar\Phi=\th_1\w\th_3+\th_2\w\th_4,\bar\Psi=\th_1\w\th_4-\th_2\w\th_3$)
$$  \n(\th_1\w\th_2)=\frac12\{\Phi(\om^4_1+\om^3_2)+\Psi(\om^1_3+\om^4_2)+\bar\Phi (-\om^4_1+\om^3_2) +\bar\Psi(-\om^1_3+\om^4_2)\}.$$
Analogously
$$  \n(\th_3\w\th_4)=\frac12\{\Phi(\om^4_1+\om^3_2)+\Psi(\om^1_3+\om^4_2)-\bar\Phi (-\om^4_1+\om^3_2) -\bar\Psi(-\om^1_3+\om^4_2)\}.$$
Note that
$\n \0=a\otimes\Phi+b\otimes\Psi$ and $\n \bO=0$, where with our assumptions $a=\a\th_1$ and $b=\a\th_2$.
On the other hand $a=\om^4_1+\om^3_2,b=\om^1_3+\om^4_2$ and
$$
\a\th_1=\om^4_1+\om^3_2,\ \a\th_2=\om^1_3+\om^4_2,\tag 2.9$$
Hence $\frac12\a\th_1=\om^4_1=\om^3_2,\ \frac12\a\th_2=\om^1_3=\om^4_2$.

\medskip
Hence, if the opposite structure is K\"ahler, then  $E_3\a=0, \G^3_{11}=\G^3_{22}=\G^3_{44}=\G^4_{21}=\G^4_{12}=0$ and
$-\G^3_{21}=\G^3_{12}=\G^4_{11}=\G^4_{22}=\frac12\a$.   We also have  $\G^4_{13}=-E_2\ln\a,\G^4_{23}=E_1\ln\a$.$\k$

{\bf Lemma  F.} {\it Let $(M,g,J)$ be a Hermitian  4-manifold with Hermitian Ricci tensor  and K\"ahler natural opposite structure. Let us assume that $|\n J|\ne 0$ on $M$. If $(M,g,\bJ)$ is a semi-symmetric surface foliated by Euclidean spaces, then  $E_4\ln\a=\frac12\a$.  On the other hand if  $(Mg,\bJ)$ is a QCH surface with Hermitian opposite structure $J$ such that $\bJ$ is a natural opposite structure for $J$  and $E_4\ln\a=\frac12\a$, then   $(Mg,\bJ)$ is semi-symmetric.}

\medskip
{\it Proof.}  If $(M,g,\bJ)$ is a semi-symmetric surface foliated by Euclidean spaces, then  $R(X,Y)E_3=R(X,Y)E_4=0$.   In particular $g(R(E_1,E_4)E_4,E_1)=0$.  Note that  $\n_{E_1}E_4=-\frac12\a E_1+ E_2\ln\a E_3, \n_{E_4}E_1=\G^2_{41}E_2$.
Note that   $\n_{E_4}E_4=0$  since $\Cal D=span\{E_3,E_4\}$ is a totally geodesic foliation. We also have
$$-\n_{E_4}\n_{E_1}E_4=\frac12 E_4\a E_1+\frac12\a\G^2_{41}E_2-E_4E_2\ln\a E_3$$  and $-\n_{[E_1,E_4]}E_4=-\frac{\a^2}4 E_1+\frac\a2E_2\ln\a E_3-E_2\ln\a\G^3_{34}E_3-\G^2_{41}\frac\a2E_2-\G^2_{41}E_1\ln\a E_3$.  Since $g(R(E_1,E_4)E_4,E_1)=0$, we obtain $E_4\a=\frac{\a^2}2$ or $E_4\ln\a=\frac12\a$.  Note that it follows also that $\G^4_{33}=\frac12\a$.   Note that $\th=-\a \th_4$ is a Lee form  of Hermitian structure $J$.    Now we have   $\tau-\kappa=6(|\th|^2+\delta\th)$, where $\tau$ is the scalar curvature of $(M,g)$, and $\kappa$  is the conformal curvature  of  $(M,g,J)$. Let us assume that  $E_4\ln\a=\frac12\a$.  Then   $\G^4_{33}=\frac12\a$. We have  $div \a E_4= tr\{X\r\n_X(\a E_4)\} =E_4\a-\a\G^4_{11}-\a\G^4_{22}+\a\G^3_{34}=-\a^2$.   Hence  $\delta \th=-\a^2=-|\th|^2$   and  $\tau-\kappa=0$. Hence our QCH surface is semi-symmetric  (see  [J-1]).
$\k$

\medskip
Let us assume that $\th_1=fdx,\th_2=fdy$.  Then  $\G^1_{32}=\frac\a2,\G^2_{41}=0$.  Then we have
$$\gather [E_1,E_4]=-\frac\a2 E_1+E_2\ln\a E_3,\tag2.10\\
[E_2,E_4]=-\frac\a2E_2-E_1\ln\a E_3,\\
[E_1,E_3]=-E_2\ln\a E_4,\\
[E_2,E_3]=E_1\ln\a E_4,\\
[E_3,E_4]=-(-E_4\ln\a+\a)E_3,\\
[E_1,E_2]=\G^1_{12}E_1-\G^2_{21}E_2+\a E_3.\endgather$$

We also have
$$\gather  d\th_1=\G^2_{11}\th_1\w\th_2+\frac12\a\th_1\w\th_4,d\th_2=-\G^1_{22}\th_1\w\th_2+\frac12\a\th_2\w\th_4,\tag2.11\\
d\th_3=-\a\th_1\w\th_2-E_2\ln\a\th_1\w\th_4+E_1\ln\a\th_2\w\th_4+(-E_4\ln\a+\a)\th_3\w\th_4,\\
d\th_4=E_2\ln\a\th_1\w\th_3-E_1\ln\a\th_2\w\th_3.\endgather$$

\medskip
 {\bf Lemma H.} {\it  The two almost Hermitian structures $J,\bJ$  on $M$,  such that the $(1,0)$  distribution of $J$ is given by  $\th_1+i\th_2=0,\th_3+i\th_4=0$ and the $(1,0)$ distribution of $\bJ$  by $\th_1+i\th_2=0,\th_3-i\th_4=0$, are integrable  if  equations (2.11) are satisfied.  The form  $d\th=-d(\a\th_4)$ is self-dual with respect to the orientation given by $\bJ$.}

 {\it Proof.}  Note that $\th_1+i\th_2=0$ if and only if  $dx+idy=0$.  On the other hand
$$\gather  d\th_3+id\th_4=-\a\th_1\w\th_2-E_2\ln\a\th_1\w\th_4+E_1\ln\a\th_2\w\th_4+(-E_4\ln \a+\a)\th_3\w\th_4+\\
iE_2\ln\a\th_1\w\th_3-iE_1\ln\a\th_2\w\th_3=-\a(\th_1+i\th_2)\w\th_2+E_2\ln\a\th_1\w(-\th_4+i\th_3)\\+E_1\ln\a\th_2\w(\th_4-i\th_3)
+(-E_4\ln \a+\a)\th_3\w\th_4=-\a(\th_1+i\th_2)\w\th_2+\\iE_2\ln\a\th_1\w(\th_3+i\th_4)-iE_1\ln\a\th_2\w(\th_3+i\th_4)
+(-E_4\ln\a+\a)(\th_3+i\th_4)\w\th_4.\endgather$$

Analogously we obtain

$$\gather  d\th_3-id\th_4=-\a\th_1\w\th_2-E_2\ln\a\th_1\w\th_4+E_1\ln\a\th_2\w\th_4+(-E_4\ln \a+\a)\th_3\w\th_4-\\
iE_2\ln\a\th_1\w\th_3+iE_1\ln\a\th_2\w\th_3=-\a(\th_1+i\th_2)\w\th_2+E_2\ln\a\th_1\w(-\th_4-i\th_3)\\+E_1\ln\a\th_2\w(\th_4+i\th_3)
+(-E_4\ln \a+\a)\th_3\w\th_4=-\a(\th_1+i\th_2)\w\th_2\\-iE_2\ln\a\th_1\w(\th_3-i\th_4)+iE_1\ln\a\th_2\w(\th_3-i\th_4)
+(-E_4\ln \a+\a)(\th_3-i\th_4)\w\th_4.\endgather$$

Note that  $d(\a\th_4)= d\a\w\th_4+\a d\th_4=E_2\a\bar{\Phi}+E_1\a\bar{\Psi}$, where $\bar{\Phi}=\th_1\w\th_3+\th_2\w\th_4,\bar{\Psi}=\th_1\w\th_4-\th_2\w\th_3$.$\k$

 Next  we give examples of QCH K\"ahler surface (see [J-1])  whose opposite almost Hermitian structure is Hermitian and not locally conformally K\"ahler.

\medskip
{\bf Theorem 1.}  {\it  Let  $U\subset \Bbb R^2$ and let $g_{\Sigma}=h^2(dx^2+dy^2)$ be a Riemannian metric on $U$, where $h:U\r \Bbb R$ is a positive function $h=h(x,y)$.  Let $\om_{\Sigma}=h^2dx\w dy$ be a volume form of $\Sigma=(U,g)$. Let $M=U\times N$, where $N=\{(z,t)\in \Bbb R^2:z<0\}$. Let us define the metric $g$  on $M$  by $g(X,Y)=z^2g_{\Sigma}(X,Y)+\th_3(X)\th_3(Y)+\th_4(X)\th_4(Y)$,   where $\th_3=-\frac z2dt+(\cos\frac12tH(x,y)+zl_2(x,y))dx+(-\sin\frac12tH(x,y)+zn_2(x,y))dy$
and
$\th_4=dz-\sin\frac12tH(x,y)dx-\cos\frac12tH(x,y)dy$ and the function $H$ satisfies the equation  $\Delta \ln H=(\ln H)_{xx}+(\ln H)_{yy}=2h^2$ on $U$, $l_2=-(\ln H)_y, n_2=(\ln H)_x$.  Then   $(M,g)$ admits a K\"ahler structure $\bJ$ with the K\"ahler form $\bO=z^2\om_{\Sigma}+\th_4\w \th_3$  and a Hermitian structure $J$ with the K\"ahler form $\Om=z^2\om_{\Sigma}+\th_3\w \th_4$.  The Ricci tensor of $(M,g)$ is $J$-invariant and $J$ is not locally conformally K\"ahler.  The Lee form of $(M,g,J)$ is $\th=-\a\th_4$,  where $\a=-\frac2z$.  The scalar curvature of $(M,g)$ is $\tau=-\frac{2\Delta\ln h}{z^2h^2}-\frac 8{z^2}$  and is equal to conformal curvature $\kappa$ of $(M,g,J)$.  $(M,g,\bO)$ is a semi-symmetric QCH K\"ahler surface.}
\medskip
{\it Proof.}
Let us take a coordinate system such that $E_1=\frac1f\p_x+k\p_z+l\p_t,E_2=\frac1f\p_y+m\p_z+n\p_t,E_3=\a\p_t,E_4=\p_z$. Then  $\th_1=fdx, \th_2=fdy, \th_4=dz-(fk)dx-(fm)dy,\th_3=\frac1\a dt-(\frac1\a lf)dx-(\frac1\a nf)dy$.  Let  $\a=-\frac 2z$. Then $E_4\ln\a=\frac12\a$ and
$$E_1\ln\a=-\frac kz,E_2\ln\a=-\frac mz.\tag2.12$$
We have $[E_1,E_4]=\frac{f_z}{f^2}\p_x-k_z\p_z-l_z\p_t=-\frac\a{2f}\p_x-\frac\a2k\p_z-\frac\a2 l\p_t-\a\frac mz\p_t$. Hence $k_z=-\frac kz,
l_z=-\frac lz-\frac{2m}{z^2}, f_z=\frac fz$. This implies   $f=zh(x,y), k=\frac{k_1(x,y,t)}z$.  On the other hand $[E_1,E_3]=-\a k_t\p_z-\a l_t\p_t+k\a_z\p_t=\frac mz\p_z$  and
$$l_t=-\frac kz, 2k_t=m.\tag 2.13$$
This yields
$$l=\frac{2m}z+\frac{l_1(x,y)}z.\tag2.14$$
Similarly $[E_2,E_3]=m\a_z\p_t-\a m_t\p_z-\a n_t\p_t=-\frac kz\p_z$
$$2m_t=-k,n_t=-\frac mz.$$
Since $[E_2,E_4]=-m_z\p_z-n_z\p_t+\frac{f_z}{f^2}\p_y=-\frac\a{2f}\p_y-\frac12\a m\p_z-\frac12\a n\p_t+\a\frac kz\p_t$, we get
$m_z=-\frac mz$.  Hence  $m=\frac{m_1(x,y,t)}z$ and
$$n=-\frac{2k}z+\frac{n_1(x,y)}z.\tag2.15$$
Let us take  $kf=\sin\frac12tH(x,y), mf=\cos\frac12t H(x,y)$,  where  $f=zh(x,y)$ and $\th_1\w\th_2=z^2h^2dx\w dy$.  Then
$lf=\frac2z\cos\frac12tH(x,y)+2l_2(x,y),nf=-\frac2z\sin\frac12tH(x,y)+2n_2(x,y)$.  Hence
$$\th_3=-\frac z2dt+(\cos\frac12tH(x,y)+zl_2(x,y))dx+(-\sin\frac12tH(x,y)+zn_2(x,y))dy$$
and
$$\th_4=dz-\sin\frac12tH(x,y)dx-\cos\frac12tH(x,y)dy.$$

Now we prove that
$d\th_3=-\a\th_1\w\th_2-E_2\ln\a\th_1\w\th_4+E_1\ln\a\th_2\w\th_4+\frac12\a\th_3\w\th_4$ and
$d\th_4=E_2\ln\a\th_1\w\th_3-E_1\ln\a\th_2\w\th_3$  if $l_2=-(\ln H)_y, n_2=(\ln H)_x$, where  $\Delta\ln H=2h^2$.

In fact,  

$-\a\th_1\w\th_2-E_2\ln\a\th_1\w\th_4+E_1\ln\a\th_2\w\th_4+\frac12\a\th_3\w\th_4=-\a f^2dx\w dy +\frac{mf}zdx\w(dz-\sin\frac12tHdx-\cos\frac12tHdy)-\frac{kf}zdy\w(dz-\sin\frac12tHdx-\cos\frac12tHdy)+\frac12\a(-\frac z2dt+(\cos\frac12t H+zl_2)dx+(-\sin \frac12 t H +zn_2)dy)\w(dz-\sin\frac12tHdx-\cos\frac12tHdy)=\frac{2f^2}zdx\w dy-\frac1z\cos^2\frac12tH^2dx\w dy +\frac1z\sin^2\frac12tH^2dy\w dx -\frac1z\sin\frac12tHdy\w dz+\frac1z\cos\frac12t Hdx\w dz-\frac 1z(-\frac z2dt\w dz+\frac z2\sin\frac12tHdt\w dx+\frac z2\cos\frac12tHdt\w dy+\cos\frac12tHdx\w dz+zl_2dx\w dz-\cos^2\frac12tH^2dx\w dy-zl_2\cos\frac12tHdx\w dy-\sin\frac12tHdy\w dz+zn_2dy\w dz+\sin^2\frac12tH^2dy\w dx -zn_2\sin\frac12tHdy\w dx)=  \frac{2f^2}zdx\w dy+\frac12dt\w dz-\frac12\sin\frac12tHdt\w dx -\frac12\cos\frac12tHdt\w dy-l_2dx\w dz+l_2 \cos\frac12tHdx\w dy-n_2dy\w dz +n_2\sin\frac12tHdy\w dx.  $

On the other hand
$$\gather d\th_3= -\frac12dz\w dt -\frac12\sin\frac12tHdt\w dx+\cos\frac12tH_ydy\w dx+l_2dz\w dx\\+zl_{2y}dy\w dx -\frac12\cos\frac12t Hdt\w dy-\sin\frac12 tH_xdx\w dy+n_2dz\w dy+zn_{2x}dx\w dy.\endgather$$
It is clear that $\th_3$ satisfies (2.11) if $l_2=-(\ln H)_y, n_2=(\ln H)_x$ and  $\Delta\ln H=2h^2$, where $\Delta f=f_{xx}+f_{yy}$, which means that $d(l_2dx+n_2dy)=2\om_{\Sigma}$.

Now  $E_2\ln\a\th_1\w\th_3-E_1\ln\a\th_2\w \th_3=-\frac{\cos\frac12tH}zdx\w(-\frac z2 dt-\sin\frac12t Hdy+zn_2dy)+\frac{\sin\frac12t H}z dy\w(-\frac z2dt+\cos\frac12t Hdx+zl_2dx)=\frac12\cos\frac12tHdx\w dt-\cos\frac12tHn_2dx\w dy-\frac12\sin\frac12tHdy\w dt+\sin\frac12tHl_2dy\w dx$.

On the other hand

$d\th_4=-\frac12\cos\frac12tHdt\w dx-\sin\frac12tH_ydy\w dx+\frac12\sin\frac12tHdt\w dy-\cos\frac12tH_xdx\w dy$.

It is clear that $\th_4$ satisfies (2.11) if $l_2=-(\ln H)_y, n_2=(\ln H)_x$.

Now we have   $d\th_1\w\th_2=d(z^2h^2dx\w dy)=2zdzh^2dx\w dy=\frac 2zdz\w\th_1\w\th_2=-\a\th_4\w\th_1\w\th_2$. From (2.11) it is also clear that  $d(\th_3\w\th_4)=-\a\th_4\w\th_1\w\th_2$.   Hence  $d\bO=0$, where $\bO(X,Y)=g(\bJ X,Y)$, and $\bO=\th_1\w\th_2-\th_3\w\th_4$. Hence in view of Lemma H $(M,g,\bJ)$ is a K\"ahler surface and the Lee form of $(M,g,J)$ is $\th=-\a\th_4$.

Now we show that $\p_t$ is a real holomorphic vector field.  Note that the opposite K\"ahler structure $\bJ$ satisfies  $\bJ\p_z=\a\p_t,\bJ\p_t=-\frac1\a\p_z,  \bJ\p_x=\p_y+fm\p_z+fn\p_t-\a fk\p_t+\frac1\a fl\p_z, \bJ\p_y=-\p_x-fk\p_z-fl\p_t-\a fm\p_t+\frac1\a fn\p_z$.
It is clear that $(L_{\p_t}\bJ)\p_z=(L_{\p_t}\bJ)\p_t=0$.  We also have  $(L_{\p_t}\bJ)\p_x=fm_t\p_z+fn_t\p_t-\a fk_t\p_t+\frac1\a fl_t\p_z=0$ and similarly    $(L_{\p_t}\bJ)\p_y=0$.  Thus   $L_{\p_t}\bJ=0$ and  $X=\p_t+i\frac1\a\p_z$ is a holomorphic vector field such that $(dt-i\a dz)(X)=2$.  The form  $\psi=\frac12(dt-i\a dz)$  satisfies  $\psi(X)=1$ and  $\psi(\p_t-i\frac1\a\p_z)=0$.   If   $\Psi$ is a holomorphic $(1,0)$ form such that  $\Psi(X)=1$,   then  $\Psi-\psi=0 mod \{dx,dy\}$.  If  $\bO$ is a K\"ahler form of K\"ahler manifold $(M,\bJ,g)$, then
$\frac12\bO^2=z^2h^2\frac1\a dx\w dy\w dz\w dt$.  We also have  $\frac{i^2}4(dx+idy)\w\overline{(dx+idy)}\w\Psi\w \overline{\Psi}=\frac{\a}4 dx\w dy\w dz\w dt$.   Hence the Ricci form  of  $(M,g,J)$  is   $\rho=-\frac12dd^c\ln \frac{z^2h^2\frac1\a dx\w dy\w dz\w dt}{\frac{\a}4 dx\w dy\w dz\w dt}=
-\frac12dd^c\ln h^2-\frac12dd^c\ln z^4=-dd^c\ln h-2dd^c\ln z=-d\bJ d\ln h-2d\bJ d\ln z= -d((\ln h)_xdy-(\ln h)_ydx)-2d(-(\ln H)_ydx+\ln H_x dy)=-\Delta\ln h dx\w dy-2\Delta(\ln H)dx\w dy$.    Consequently $E_1,E_2,E_3,E_4$ are eigenfields of the Ricci tensor and  $E_3,E_4$ corresponds to the eigenvalue $0$.
Since  $\rho=-(\Delta\ln h+4h^2) dx\w dy=-\frac{\Delta\ln h+4h^2}{z^2h^2}z^2h^2 dx\w dy$, it follows that $\frac12\tau=-\frac{\Delta\ln h}{z^2h^2}-\frac4{z^2}$.$\k$

Now we give examples of non-semi-symmetric QCH surfaces with non l.c.K. opposite Hermitian structure.

{\bf Theorem  2.}  {\it  Let  $U\subset \Bbb R^2$ and let $g_{\Sigma}=h^2(dx^2+dy^2)$ be a Riemannian metric on $U$, where $h:U\r \Bbb R$ is a positive function $h=h(x,y)$.  Let $\om_{\Sigma}=h^2dx\w dy$ be a volume form of $\Sigma=(U,g)$. Let $M=U\times N$, where $N=\{(z,t)\in \Bbb R^2:0<z<\pi\}$. Let us define the metric $g$  on $M$  by $g(X,Y)=(\cos\frac12 z)^2g_{\Sigma}(X,Y)+\th_3(X)\th_3(Y)+\th_4(X)\th_4(Y)$,   where $ \th_3=\sin zdt-(\cos t\cos zH(x,y)+\sin zl_2(x,y))dx-(-\sin t\cos zH(x,y)+\sin zn_2(x,y))dy$,
$\th_4=dz-\sin tH(x,y)dx-\cos tH(x,y)dy$ and the function $H$ satisfies the equation  $\Delta \ln H=(\ln H)_{xx}+(\ln H)_{yy}=\frac12h^2-H^2$ on $U$, $l_2=-(\ln H)_y, n_2=(\ln H)_x$.  Then   $(M,g)$ admits a K\"ahler structure $\bJ$ with the K\"ahler form $\bO=(\cos\frac12z)^2\om_{\Sigma}+\th_4\w \th_3$  and a Hermitian structure $J$ with the K\"ahler form $\Om=(\cos\frac12z)^2\om_{\Sigma}+\th_3\w \th_4$.  The Ricci tensor of $(M,g)$ is $J$-invariant and $J$ is not locally conformally K\"ahler.  The Lee form of $(M,g,J)$ is $\th=-\a\th_4$,  where $\a=\tan\frac z2$.  }

\medskip

{\it Proof.} Let us take a coordinate system such that $E_1=\frac1f\p_x+k\p_z+l\p_t,E_2=\frac1f\p_y+m\p_z+n\p_t,E_3=\frac1\beta\p_t,E_4=\p_z$. Then  $\th_1=fdx, \th_2=fdy, \th_4=dz-(fk)dx-(fm)dy,\th_3=\beta dt-(\beta lf)dx-(\beta nf)dy$.  Let  $\a=\tan\frac z2,\beta=\sin z$. Then

$$E_1\ln\a=\frac k{\sin z},E_2\ln\a=\frac m{\sin z}.\tag2.16$$
We have $[E_1,E_4]=\frac{f_z}{f^2}\p_x-k_z\p_z-l_z\p_t=-\frac\a{2f}\p_x-\frac\a2k\p_z-\frac\a2 l\p_t+\frac m{\sin^2z}\p_t$. Hence $k_z=\frac12\tan\frac12z k,
l_z=\frac12\tan\frac12zl-\frac{m}{sin^2z}, f_z=-\frac12\tan\frac12zf$. This implies   $$f=\cos\frac12zh(x,y), k=\frac{k_1(x,y,t)}{\cos\frac12z}.$$ On the other hand $[E_1,E_3]=-\frac1\beta( k_t\p_z+ l_t\p_t)-k\frac{\beta_z}{\beta^2}\p_t=-\frac m{\sin z}\p_z$  and
$$l_t= -k\cot z, k_t=m.\tag 2.17$$
This yields
$$l=m\cot z+\frac{l_1(x,y)}{\cos\frac12 z}.\tag2.18$$
Similarly $[E_2,E_3]=-m\frac{\beta_z}{\beta^2}\p_t-\frac1{\beta} (m_t\p_z+ n_t\p_t)=-\frac k{\sin z}\p_z$,
$$m_t=-k,n_t= -m\cot z.$$
Since $[E_2,E_4]=-m_z\p_z-n_z\p_y+\frac{f_z}{f^2}\p_t=-\frac\a{2f}\p_y-\frac12\a m\p_z-\frac12\a n\p_t-\frac k{\sin^2z}\p_t$, we get
$m_z=\frac12\tan\frac12z m$.  Hence  $m=\frac{m_1(x,y,t)}{\cos\frac12z}$ and
$$n=-k\cot z+\frac{n_1(x,y)}{\cos\frac12z}.\tag2.19$$
Let us take  $kf=\sin tH(x,y), mf=\cos t H(x,y)$,  where  $f=\cos\frac12zh(x,y)$ and $\th_1\w\th_2=\cos\frac12z^2h^2dx\w dy$.  Then
$$lf=\cot z\cos tH(x,y)+l_2(x,y),nf=-\cot z\sin tH(x,y)+n_2(x,y).$$  Hence
$$\gather \th_3=\sin zdt-(\cos t\cos zH(x,y)+\sin zl_2(x,y))dx-(-\sin t\cos zH(x,y)\\+\sin zn_2(x,y))dy\endgather$$
and
$$\th_4=dz-\sin tH(x,y)dx-\cos tH(x,y)dy.$$

Now we prove that
$d\th_3=-\a\th_1\w\th_2-E_2\ln\a\th_1\w\th_4+E_1\ln\a\th_2\w\th_4+(-E_4\ln\a+\a)\th_3\w\th_4$ and
$d\th_4=E_2\ln\a\th_1\w\th_3-E_1\ln\a\th_2\w\th_3$  if $l_2=-(\ln H)_y, n_2=(\ln H)_x$, where  $\Delta\ln H=-H^2+\frac12h^2$.

In fact,  $$\gather-\a\th_1\w\th_2-E_2\ln\a\th_1\w\th_4+E_1\ln\a\th_2\w\th_4+(-E_4\ln\a+\a)\th_3\w\th_4=\\-\tan \frac12z f^2dx\w dy -\frac{mf}{\sin z}dx\w(dz-\sin tHdx-\cos tHdy)\\+\frac{kf}{\sin z}dy\w(dz-\sin tHdx-\cos tHdy)-\frac{\cos z}{\sin z}(\sin zdt-\cos z \cos tH dx\\-\sin zl_2 dx+\sin  t\cos z H dy -\sin zn_2dy)\w(dz-\sin tHdx-\cos tHdy)\\=-\frac12\sin zh^2dx\w dy+\frac1{\sin z}\cos^2 tH^2dx\w dy -\frac1{\sin z}\sin^2tH^2dy\w dx\\ +\frac{\sin tH}{\sin z}dy\w dz-\frac{\cos tH}{\sin z}dx\w dz-\frac {\cos z}{\sin z}(\sin zdt\w dz\\-\sin z\sin tHdt\w dx-\sin z\cos tHdt\w dy\\-\cos z\cos tHdx\w dz-\sin zl_2dx\w dz+\cos z\cos^2tH^2dx\w dy+\\ \sin zl_2\cos tHdx\w dy+\cos z\sin tHdy\w dz\\-\sin zn_2dy\w dz-\cos z\sin^2tH^2dy\w dx +n_2\sin z\sin tHdy\w dx). \endgather$$

On the other hand
$$\gather d\th_3= \cos zdz\w dt +\sin z\cos tHdz\w dx+\cos z\sin tHdt\w dx-\\\cos z\cos tH_ydy\w dx-l_2\cos zdz\w dx-\sin zl_{2y}dy\w dx +\cos z\cos t Hdt\w dy\\-\sin t\sin zHdz\w dy+\sin t\cos zH_xdx\w dy-\cos zn_2dz\w dy-\sin zn_{2x}dx\w dy.\endgather$$

It is clear that $\th_3$ satisfies (2.11) if $l_2=-(\ln H)_y, n_2=(\ln H)_x$ and  $\Delta\ln H=\frac12h^2-H^2$, where $\Delta f=f_{xx}+f_{yy}$.

Now  $E_2\ln\a\th_1\w\th_3-E_1\ln\a\th_2\w \th_3=\frac{\cos tH}{\sin z}dx\w(\sin z dt+\sin t\cos z Hdy-\sin zn_2dy)-\frac{\sin t H}{\sin z} dy\w(\sin zdt-\cos z\cos t Hdx-\sin zl_2dx)=cos tHdx\w dt-\cos tHn_2dx\w dy +\frac{\cos t\sin t\cos zH^2}{\sin z}dx\w dy-\sin t Hdy\w dt+\frac{\sin t \cos t\cos zH^2}{\sin z}dy\w dx+\sin tHl_2dy\w dx.$.

On the other hand

$d\th_4=-\cos tHdt\w dx-\sin tH_ydy\w dx+\sin tHdt\w dy-\cos tH_xdx\w dy$.

It is clear that $\th_4$ satisfies (2.11) if $l_2=-(\ln H)_y, n_2=(\ln H)_x$.

Now we have   $d\th_1\w\th_2=d(\cos\frac12z^2h^2dx\w dy)=-\sin\frac12z\cos\frac12zdzh^2dx\w dy=-\tan\frac12zdz\w\th_1\w\th_2=-\a\th_4\w\th_1\w\th_2$. From (2.11) it is also clear that  $d\th_3\w\th_4=-\a\th_4\w\th_1\w\th_2$.   Hence  $d\bO=0$, where $\bO(X,Y)=g(\bJ X,Y)$ and $\bO=\th_1\w\th_2-\th_3\w\th_4$. Hence in view of Lemma H $(M,g,\bJ)$ is a K\"ahler surface and the Lee form of $(M,g,J)$ is $\th=-\a\th_4$.

Now we show that $\De=span\{E_3,E_4\}$ is a conformal foliation.  Note that the opposite K\"ahler structure $\bJ$ satisfies  $\bJ\p_z=\frac1\beta\p_t,\bJ\p_t=-\beta\p_z,  \bJ\p_x=\p_y+fm\p_z+fn\p_t-\frac1\beta fk\p_t+\beta fl\p_z, \bJ\p_y=-\p_x-fk\p_z-fl\p_t-\frac1\beta fm\p_t+\beta fn\p_z$.  Hence  $(L_{\p_t}\bJ)\p_x,(L_{\p_z}\bJ)\p_x,(L_{\p_t}\bJ)\p_y,(L_{\p_z}\bJ)\p_y\in\De$,  which means that $L_{\xi}\bJ(\DE)\subset \De$ for $\xi\in\G(\De)$. Thus foliation $\De$ is conformal and since  $d\th^-=0$, it follows from [J-3] that $(M,g,\bJ)$ is a QCH K\"ahler surface.$\k$
\bigskip
{\bf Theorem  3.}  {\it  Let  $U\subset \Bbb R^2$ and let $g_{\Sigma}=h^2(dx^2+dy^2)$ be a Riemannian metric on $U$, where $h:U\r \Bbb R$ is a positive function $h=h(x,y)$.  Let $\om_{\Sigma}=h^2dx\w dy$ be a volume form of $\Sigma=(U,g)$. Let $M=U\times N$, where $N=\{(z,t)\in \Bbb R^2:z<0\}$. Let us define the metric $g$  on $M$  by $g(X,Y)=(\sinh\frac12 z)^2g_{\Sigma}(X,Y)+\th_3(X)\th_3(Y)+\th_4(X)\th_4(Y)$,   where $$\gather \th_3=\sinh zdt-(-\sin t\cosh z H(x,y)+\sinh z l_2(x,y))dx - \\(\cos t \cosh z H(x,y)+\sinh zn_2(x,y))dy\endgather$$and
  $\th_4=dz-\cos tH(x,y)dx-\sin tH(x,y)dy$ and the function $H$ satisfies the equation  $\Delta \ln H=(\ln H)_{xx}+(\ln H)_{yy}=(\frac12h^2+H^2)$ on $U$, $l_2=(\ln H)_y, n_2=-(\ln H)_x$.  Then   $(M,g)$ admits a K\"ahler structure $\bJ$ with the K\"ahler form $$\bO = (\sinh\frac12z)^2\om_{\Sigma} +\th_4\w \th_3$$  and a Hermitian structure $J$ with the K\"ahler form $\Om=(\sinh\frac12z)^2\om_{\Sigma}+ \th_3\w \th_4$.  The Ricci tensor of $(M,g)$ is $J$-invariant and $J$ is not locally conformally K\"ahler.  The Lee form of $(M,g,J)$ is $\th=-\a\th_4$,  where $\a=-\coth\frac z2$.  }

\medskip

{\it Proof.} Let us now take a coordinate system such that $E_1=\frac1f\p_x+k\p_z+l\p_t,E_2=\frac1f\p_y+m\p_z+n\p_t,E_3=\frac1\beta\p_t,E_4=\p_z$. Then  $\th_1=fdx, \th_2=fdy, \th_4=dz-(fk)dx-(fm)dy,\th_3=\beta dt-(\beta lf)dx-(\beta nf)dy$.  Let  $\a=-\coth\frac z2,\beta=\sinh z$. Then

$$E_1\ln\a=-\frac k{\sinh z},E_2\ln\a=-\frac m{\sinh z}.\tag2.20$$
We have $[E_1,E_4]=\frac{f_z}{f^2}\p_x-k_z\p_z-l_z\p_t=-\frac\a{2f}\p_x-\frac\a2k\p_z-\frac\a2 l\p_t-\frac m{\sinh^2z}\p_t$. Hence $k_z=-\frac12\coth\frac12z k,
l_z=-\frac12\coth\frac12zl+\frac{m}{\sinh^2z}, f_z=\frac12\coth\frac12zf$.  This implies   $f=\sinh\frac12zh(x,y), k=\frac{k_1(x,y,t)}{\sinh\frac12z}$.  On the other hand $[E_1,E_3]=-\frac1\beta( k_t\p_z+ l_t\p_t)-k\frac{\beta_z}{\beta^2}\p_t=\frac m{\sinh z}\p_z$  and
$$l_t= -k\coth z, k_t=-m.\tag 2.21$$
This yields
$$l=-m\coth z+\frac{l_1(x,y)}{\sinh\frac12 z}.\tag2.22$$

Similarly $[E_2,E_3]=-m\frac{\beta_z}{\beta^2}\p_t-\frac1{\beta} (m_t\p_z+ n_t\p_t)=-\frac k{\sinh z}\p_z$
$$m_t=k,n_t= -m\coth z.$$
Since $[E_2,E_4]=-m_z\p_z-n_z\p_t+\frac{f_z}{f^2}\p_y=-\frac\a{2f}\p_y-\frac12\a m\p_z-\frac12\a n\p_t+\frac k{\sinh^2z}\p_t$, we get
$m_z=-\frac12\coth\frac12z m$.  Hence  $m=\frac{m_1(x,y,t)}{\sinh\frac12z}$ and
$$n=k\coth z+\frac{n_1(x,y)}{\sinh\frac12z}.\tag2.23$$
Let us take  $kf=\cos tH(x,y), mf=\sin t H(x,y)$,  where  $f=\sinh\frac12zh(x,y)$ and $\th_1\w\th_2=\sinh\frac12z^2h^2dx\w dy$.  Then
$$lf\beta=-\sin t\cosh zH+\sinh zl_2(x,y),nf\beta=\cosh z\cos tH(x,y)+\sinh zn_2(x,y).$$  Hence
$$\gather \th_3=\sinh zdt-(-\sin t\cosh zH(x,y)+\sinh zl_2(x,y))dx-(\cos t\cosh zH(x,y)\\+\sinh zn_2(x,y))dy\endgather$$
and
$$\th_4=dz-\cos tH(x,y)dx-\sin tH(x,y)dy.$$

Now we prove that
$d\th_3=-\a\th_1\w\th_2-E_2\ln\a\th_1\w\th_4+E_1\ln\a\th_2\w\th_4+(-E_4\ln\a+\a)\th_3\w\th_4$ and
$d\th_4=E_2\ln\a\th_1\w\th_3-E_1\ln\a\th_2\w\th_3$  if $l_2=-(\ln H)_y, n_2=(\ln H)_x$, where  $\Delta\ln H=-(H^2+\frac12h^2)$.

In fact,  $$\gather-\a\th_1\w\th_2-E_2\ln\a\th_1\w\th_4+E_1\ln\a\th_2\w\th_4+(-E_4\ln\a+\a)\th_3\w\th_4=\\\coth \frac12z f^2dx\w dy +\frac{mf}{\sinh z}dx\w(dz-\sin tHdy)\\-\frac{kf}{\sinh z}dy\w(dz-\cos tHdx)-\frac{\cosh z}{\sinh z}(\sinh zdt+\cosh z \sin tH dx\\-\sinh zl_2 dx-\cos  t\cosh z H dy -\sinh zn_2)dy)\w(dz-\cos tHdx-\sin tHdy)\\=\frac12\sinh zh^2dx\w dy-\frac1{\sinh z}\sin^2 tH^2dx\w dy +\frac1{\sinh z}\cos^2tH^2dy\w dx\\ -\frac{\cos tH}{\sinh z}dy\w dz+\frac{\sin tH}{\sinh z}dx\w dz-\frac {\cosh z}{\sinh z}(\sinh zdt\w dz\\-\sinh z\cos tHdt\w dx-\sinh z\sin tHdt\w dy\\+\cosh z\sin tHdx\w dz-\sinh zl_2dx\w dz-\cosh z\sin^2tH^2dx\w dy+\\ \sinh zl_2\sin tHdx\w dy-\cosh z\cos tHdy\w dz\\-\sinh zn_2dy\w dz+\cosh z\cos^2tH^2dy\w dx +n_2\sinh z\cos tHdy\w dx). \endgather$$

On the other hand
$$\gather d\th_3= \cosh zdz\w dt +\sinh z\sin tHdz\w dx+\cosh z\cos tHdt\w dx+\\ \cosh z\sin tH_ydy\w dx-l_2\cosh zdz\w dx-\sinh zl_{2y}dy\w dx +\cosh z\sin t Hdt\w dy\\-\cos t\sinh zHdz\w dy-\cos t\cosh zH_xdx\w dy-\cosh zn_2dz\w dy\\-\sinh zn_{2x}dx\w dy.\endgather$$

It is clear that $\th_3$ satisfies (2.11) if $l_2=(\ln H)_y, n_2=-(\ln H)_x$ and  $\Delta\ln H=(\frac12h^2+H^2)$, where $\Delta f=f_{xx}+f_{yy}$.

Now  $E_2\ln\a\th_1\w\th_3-E_1\ln\a\th_2\w \th_3=-\frac{\sin tH}{\sinh z}dx\w(\sinh z dt-\cos t\cosh z Hdy-\sinh zn_2dy)+\frac{\cos t H}{\sinh z} dy\w(\sinh zdt+\cosh z\sin t Hdx-\sinh zl_2dx)=-sin tHdx\w dt+\sin tHn_2dx\w dy +\frac{\cos t\sin t\cosh zH^2}{\sinh z}dx\w dy+\cos t Hdy\w dt+\frac{\sin t \cos t\cosh zH^2}{\sinh z}dy\w dx-\cos tHl_2dy\w dx=-sin tHdx\w dt+\sin tHn_2dx\w dy+\cos t Hdy\w dt-\cos tHl_2dy\w dx$.

On the other hand

$d\th_4=\sin tHdt\w dx-\cos tH_ydy\w dx-\cos tHdt\w dy-\sin tH_xdx\w dy$.

It is clear that $\th_4$ satisfies (2.11) if $l_2=(\ln H)_y, n_2=-(\ln H)_x$.

Now we have   $d\th_1\w\th_2=d(\sinh\frac12z^2h^2dx\w dy)=\sinh\frac12z\cosh\frac12zdzh^2dx\w dy=-(-\coth\frac12z)dz\w\th_1\w\th_2=-\a\th_4\w\th_1\w\th_2$. From (2.11) it is also clear that  $d\th_3\w\th_4=-\a\th_4\w\th_1\w\th_2$.   Hence  $d\bO=0$, where $\bO(X,Y)=g(\bJ X,Y)$ and $\bO=\th_1\w\th_2-\th_3\w\th_4$. Hence in view of Lemma H $(M,g,\bJ)$ is a K\"ahler surface and the Lee form of $(M,g,J)$ is $\th=-\a\th_4$.  As above one can show  that  $(M,g,\bJ)$  is a QCH K\"ahler surface.$\k$

\par
\bigskip

\centerline{\bf References.}

\par
\medskip
\cite{A-G-1} V. Apostolov and P. Gauduchon, {\it The Riemannian Goldberg-Sachs theorem}, Int. J. Math. {\bf 8}, No.4,(1997), 421-439.
\par
\par
\medskip
\cite{D} T.C. Draghici, {\it Almost K\"ahler 4-manifolds with J-invariant Ricci tensor}, Houston J. Math. {\bf 25}, No. 1,(1999), 133-145.
\par
\medskip
\cite{G-H} A. Gray and L.M. Harvella, {\it The sixteen classes of almost Hermitian manifolds} Ann. Mat. Pura    Appl.{\bf 123}, (1980), 35-58.
\par
\medskip
\cite{J-1} W. Jelonek {\it K\"ahler surfaces with quasi constant holomorphic curvature} Glasgow Math.J.58,(2016), 503-512.
\par
\medskip
\cite{J-2} W. Jelonek {\it Semi-symmetric  K\"ahler surfaces } Colloq. Math.J.148,(2017), 1-12.
\par
\medskip
\cite{J-3} W. Jelonek {\it Complex foliations and  K\"ahler QCH surfaces }  (arxiv: 1412.3363 v4 [math DG] to appear in Colloq. Math.)
\par
\medskip
\cite{J-4} W. Jelonek {\it Einstein Hermitian and anti-Hermitian 4-manifolds } Ann.   \newline Polon. Math. 81.1 (2003), 7-24.
\par
\medskip
\cite{M} O. Muskarov  {\it On hermitian surfaces with J-invariant Ricci tensor} Journal of Geom. 72,(2001),151-156.

\par
\medskip

 Institute of Mathematics

 Cracow  University of Technology

 Warszawska 24

 31-155 Krak\'ow,POLAND.

\enddocument